\documentclass[lettersize,journal]{IEEEtran}
\usepackage{amsmath,amsfonts}
\usepackage{algorithmic}
\usepackage{algorithm}
\usepackage{array}
\usepackage[caption=false,font=normalsize,labelfont=sf,textfont=sf]{subfig}
\usepackage{textcomp}
\usepackage{stfloats}
\usepackage{url}
\usepackage{verbatim}
\usepackage{graphicx}
\usepackage{cite}
\begin{document}

\title{A High-order Nyström-based Scheme Explicitly Enforcing Surface Density Continuity for the Electric Field Integral Equation }

\author{Jin Hu, \IEEEmembership{Student Member, IEEE} and Constantine Sideris, \IEEEmembership{Senior Member, IEEE}
\thanks{The authors gratefully acknowledge support by the Air Force Office of Scientific Research (FA9550-20-1-0087) and the National Science Foundation (CCF-2047433).}}



\maketitle

\begin{abstract}
This paper introduces an efficient approach for solving the Electric Field Integral Equation (EFIE) with high-order accuracy by explicitly enforcing the continuity of the impressed current densities across boundaries of the surface patch discretization. The integral operator involved is discretized via
a Nyström-collocation approach based on Chebyshev polynomial expansion within each patch and a closed quadrature rule is utilized such that the discretization points inside one patch coincide with those inside another patch on the shared boundary of those two patches. The continuity enforcement is achieved by constructing a mapping from those coninciding points to a vector containing unique discretization points used in the GMRES iterative solver. The proposed approach is applied to the scattering of several different geometries including a sphere, a cube, a NURBS model imported from CAD  software, and a dipole structure and results are compared with the Magnetic Field Integral
Equation (MFIE) and the EFIE without enforcing continuity to illustrate the effectiveness of the approach.
\end{abstract}


\section{Introduction}
Boundary Integral Equation (BIE) methods are among the most popular numerical approaches for solving Maxwell's equations. They are used extensively for radiation and scattering analysis in many applications, including antennas, radio-frequency (RF) circuits and nanophotonics devices and are attractive over competing volumetric methods such as finite-difference and finite element methods due to only requiring surface mesh discretizations. 
Most solution approaches for BIEs leverage low-order Method-of-Moment (MoM) discretizations using RWG basis functions~\cite{rao1982electromagnetic}. However, the solution accuracy may be limited both due to the low-order basis expansion of the current density as well as the flat triangular discretization used for approximating the original geometry. Further, MoM relies on Galerkin testing, which results in requiring costly 4D integrals to build the impedance matrix. Unlike MoM, the Locally Corrected Nystrom (LCN) method~\cite{canino1998} leverages point-matching for testing and can still achieve high-order accuracy when the locally corrected weights accounting for near interactions are evaluated properly. However, the calculation of these weights can be time consuming. Recently, a new Nystr\"om method was proposed for the acoustic~\cite{bruno2020chebyshev} and electromagnetic scattering problems~\cite{hu2021chebyshev, garza2021aces, garza2022acs, hu2023uniaxial}, which uses Chebyshev polynomials to approximate the current on the surface and uses a singularity cancellation approach to precompute the local corrections rapidly with high accuracy. 

Despite their advantages over MoM, an issue for Nystr\"om based methods when solving the EFIE problems is the  continuity of the current density across mesh element edges is not explicitly enforced, while the divergence-conforming RWG function used in MoM inherently impose the continuity. As a result, line integrals are usually incorporated in the EFIE formulation to compensate for the potential discontinuity for a Nystr\"om discretization  \cite{gedney2003deriving} or a set of line charges have to be introduced as extra unknowns to augment to EFIE formulations \cite{young2012high} to solve this issue. Another approach is to construct a special set of basis functions by performing singular value decomposition of constraint matrices obtained by imposing the continuity at every node of each patch boundary\cite{hendijani2015constrained}.

In this letter, a new scheme is proposed for explicitly enforcing the current density in the EFIE solved with the Chebyshev-based Nystr\"om method (CBIE) \cite{hu2021chebyshev} by choosing a set of discretization points on each patch that also locates points on patch boundaries and corners and requiring that points from adjacent patches which physically coincide take on the same density. 
No extra line integrals, line charges, or SVDs are required, and only weakly singular integrals involving the scalar free-space Green function are computed, leading to very efficient implementation. 
To showcase the effectiveness of this approach, various numerical examples are presented and compared against results obtained with an MFIE formulation, the EFIE formulation without enforcing continuity, and the commercial solver Ansys HFSS.

\section{EFIE with Explicit Continuity Enforcement\label{sec:formulation}}
Consider the scattering of the time-harmonic electromagnetic wave by a PEC object embedded in a 3D homogeneous medium characterized by permittivity $\varepsilon$ and permeability $\mu$. The incident excitation and the scattered field are denoted by $\left(\mathbf{E}^\text{inc},\mathbf{H}^\text{inc}\right)$ and $\left(\mathbf{E}^\text{scat},\mathbf{H}^\text{scat}\right)$ respectively. Based on the representation theorem, the scattered field can be expressed in terms of the induced electric current density $\mathbf{J}=\mathbf{\hat{n}}\times\mathbf{H}$ on the surface $\Gamma$ of the object:
\begin{equation}
     \mathbf{E}^{\text{scat}} = ik\eta \int_{\Gamma}(\overline{\mathbf{I}}+\frac{\nabla\nabla}{k^{2}})G(\mathbf{r},\mathbf{r}')\cdot \mathbf{J}(\mathbf{r}')d\sigma(\mathbf{r'})
\end{equation}
where $G\left(\mathbf{r},\mathbf{r}'\right)=\exp\left(ik\left|\mathbf{r}-\mathbf{r}'\right|\right)/(4\pi\left|\mathbf{r}-\mathbf{r}'\right|)$ is the scalar Green's function, $k=\omega\sqrt{\mu\varepsilon}$ is the wavenumber, $\eta=\sqrt{{\mu}/{\varepsilon}}$ is the impedace of the medium, and $\overline{\mathbf{I}}$ is the unit dyadic. The electric field integral equation (EFIE) can be derived by enforcing the PEC boundary conditions requiring the tangential electric fields to be zero:
\begin{equation}
     \mathcal{T}\mathbf{J} = -\mathbf{\hat{n}} \times \mathbf{E}^{\text{inc}}
    \label{eq:efie}
\end{equation}
where the operator $\mathcal{T}$ is defined as
\begin{equation}
    \mathcal{T}\mathbf{J} = ik\eta \mathbf{\hat{n}} \times \int_{\Gamma}(\overline{\mathbf{I}}+\frac{\nabla\nabla}{k^{2}})G(\mathbf{r},\mathbf{r}')\cdot \mathbf{J}(\mathbf{r}')d\sigma(\mathbf{r'})
\end{equation}
Note that the $\nabla\nabla G$ term in the kernel of the integral is hypersingular with $O(1/\mathrm{R}^3)$ singularity, which complicates accurate numerical evaluation. To alleviate this singularity, vector identities can be utilized to transfer one of the $\nabla$ operators to the current density $\mathbf{J}$ and factor the remaining $\nabla$ operator out of the integral such that the resulting surface integral has a weakly singular kernel $G$:
\begin{equation}
    \begin{aligned}   
    \int_{\Gamma}\nabla\nabla G\cdot \mathbf{J}d\sigma(\mathbf{r'}) = \nabla\int_{\Gamma} G \nabla'_s\cdot \mathbf{J}d\sigma(\mathbf{r'}) - \int_{C}\nabla G  \mathbf{J} \cdot \mathbf{e}dl(\mathbf{r'})
    \end{aligned}
    \label{eq:efie_weaksing}
\end{equation}
where $C$ is the bounding contour of the surface $\Gamma$ and $\mathbf{e} = \mathbf{e}(\mathbf{r'})$ is the unit normal vector to $C$ lying in the tangent plane of $\Gamma$ at $\mathbf{r'}$. Note that the line integral over $C$ in the RHS of (\ref{eq:efie_weaksing}) has kernel $\nabla G$ with $O(1/\mathrm{R}^2)$ singularity, which is strongly singular. However, if $\Gamma$ is a smooth closed surface and the component of the current density $\mathbf{J}$ normal to any boundary is continuous over the whole surface, the contribution of this term is always zero. To see this, suppose the surface $\Gamma$ is discretized into $M$ non-overlapping quadrilateral patches: $\Gamma = \cup_{i=1}^{M}{\Gamma_{i}}$ with $\Gamma_{i} \cap \Gamma_{j} =  \emptyset $ if $i\neq j$. The contributions to the line integral from the edge shared by adjacent patches $\Gamma_{m}$ and $\Gamma_{n}$ is
\begin{equation}
    \begin{aligned}   
    &&\int_{C_{mn}}\nabla G  \mathbf{J}_{m} \cdot \mathbf{e}_{mn}dl(\mathbf{r'})+\int_{C_{mn}}\nabla G  \mathbf{J}_{n} \cdot \mathbf{e}_{nm}dl(\mathbf{r'}) \\
    && =\int_{C_{mn}}\nabla G  (\mathbf{J}_{m} \cdot \mathbf{e}_{mn} - \mathbf{J}_{n} \cdot \mathbf{e}_{mn})dl(\mathbf{r'})
    \end{aligned}
\end{equation}
where $\mathbf{J}_{m}$, $\mathbf{J}_{n}$ denote the current density on $\Gamma_{m}$, $\Gamma_{n}$ respectively and $\mathbf{e}_{mn}$ points from $\Gamma_{m}$ to $\Gamma_{n}$. Thus, if the normal component of the current density is continuous across the shared edge $C_{mn}$ (i.e., $\mathbf{J}_{m}\cdot \mathbf{e}_{mn} = \mathbf{J}_{n}\cdot \mathbf{e}_{mn}$), the contributions from the patches cancel and the line integral evaluates to exactly 0. However, in order to avoid computing this line integral in the Nystr\"om discretization scheme, the discrete representation used for approximating the real physical induced current density on the surface must also satisfy the same normal continuity condition (i.e., it must be divergence free). 

\begin{figure*}
  \centering\includegraphics[width=0.8\textwidth,height=4.8cm]{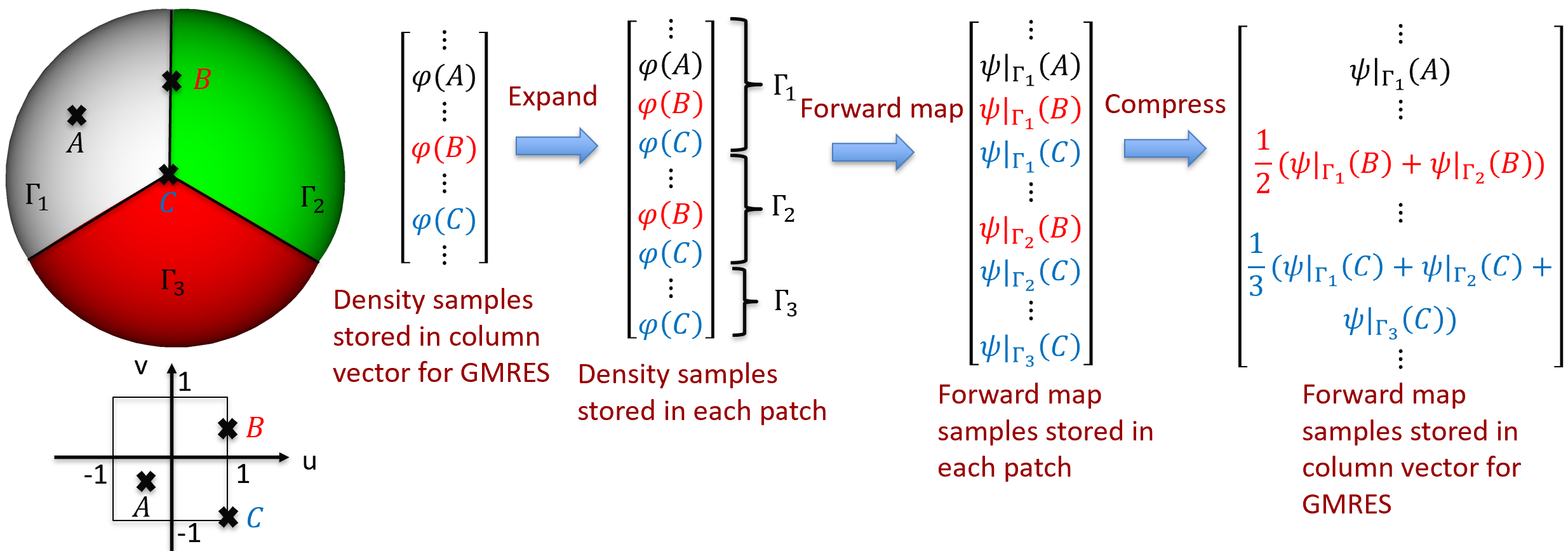}
  \caption{Illustration for the whole forward map used in GMRES iterative solver with $\varphi$ and $\phi$ representing either one of the contravariant component of $\mathbf{J}$ and $\mathcal{T}\mathbf{J}$ in equation \eqref{eq:efie}.}
  \label{fig:explain}
\end{figure*}

Unlike the Fejer's first quadrature used in \cite{hu2021chebyshev} for the MFIE, which is an open integration rule, instead we use Clenshaw–Curtis quadrature, which also includes the endpoints. 
Suppose for each one of $\Gamma_{m}$, $\Gamma_{n}$ patches mentioned above, a mapping is defined to map the square UV $[-1,1]\times[1,1]$ parametric domain to Cartesian coordinates. The Clenshaw-Curtis quadrature nodes along the U or V direction are:
\begin{equation}
    \begin{aligned}   
    x_i = \cos\left(\pi\frac{i}{N-1}\right),\quad i=0,...,N-1, 
    \end{aligned}
\end{equation}
As with the original Chebyshev method, the unknowns over each patch are located on the same grid as the quadrature rule. The current density on $\Gamma_{m}$ can be thus be discretized on these quadrature nodes: \begin{multline}
    \mathcal{J}^m = [ 
    J^{m,u}(u_0,v_0),
    \dots ,
    J^{m,u}(u_{N-1},v_{N-1}), \\ 
    J^{m,v}(u_0,v_0), 
    \dots ,
    J^{m,v}(u_{N-1},v_{N-1})]^T
    \label{eq:Jp}
\end{multline} 
where $J^{m,u}$ and $J^{m,v}$ denote the contravariant components of $\mathbf{J}_{m}$ respectively, $u_l$ and $v_k$ are the discretization points on the Chebyshev grid corresponding to the $x_i$ nodes: $\left.u_l = x_l\right|l=0,\dots,N-1,\left.v_k = x_k\right|k=0,\dots,N-1$ for an $N\times N$ discretization. Note that unlike the method in~\cite{hu2021chebyshev}, since the unknowns are located on the Clenshaw-Curtis quadrature nodes which include the endpoints, each patch also has unknowns on its boundaries and corners. This implies that at shared edges or corners between different patches, points belonging to different patches may be located in the same physical locations. Therefore, the discretization points within each patch can be classified into three types as shown in Fig.~\ref{fig:explain}: Type A: $u_l \neq \pm1$ and $v_k \neq \pm1$ corresponds to points inside the patch that do not coincide with any points from other patches (no duplicates). Type B: $u_l = \pm1$ or $v_k = \pm1$ corresponds to points on the edges that share the same location with exactly one point from an adjacent patch (multiplicity: 2). Type C: $u_l = \pm1$ and $v_k = \pm1$ corresponds to points on corners that share the same location with a point from two or three other patches (multiplicity: 2 or 3). 

In order to arrive at a uniquely solvable system, the duplicated points (type B and C) must be taken care of. Fig.~\ref{fig:explain} illustrates visually the whole matrix vector multiplication procedure of the proposed continuity-enforced approach. Please refer to~\cite{hu2021chebyshev} for details on the Chebyshev-based discretization method for the source-to-target point forward mapping approach. We define a simple tall 0-1 projection matrix which maps elements from the unknown vector in a 1-to-$M$ fashion for any type B and C points, where $M$ is the number of patches sharing the same point and 1-to-1 for all of the remaining interior (type A) points. In order to enforce a continuous density, points shared by multiple patches must be assigned to the same amplitude, and this can be done by having multiple nonzero entries in columns of the projection matrix which map an element of the unique current density vector to multiple type B or C points (e.g., in Fig.~\ref{fig:explain}, $\varphi(B)$ is copied into patches $\Gamma_{1}$ and $\Gamma_{2}$ and $\varphi(C)$ is copied into patches $\Gamma_{1}$, $\Gamma_{2}$ and $\Gamma_{3}$). After the BIE operator is applied, which maps the current density amplitudes on each source point back to scattered fields on the same set of target points, the resulting output vector must be compressed back down to the same size as the original input vector of unique unknowns in order to yield a full rank square system. This is done by defining an output compression matrix which consists of each row of transpose of the input projection matrix multiplied by a scalar factor such that the elements in each row of the compression matrix sum to 1. Thus, this output matrix effectively takes the average of target point evaluations which have been duplicated due to being at the same physical location (e.g., $\psi|_{\Gamma_{1}}(B)$ and $\psi|_{\Gamma_{2}}(B)$) or the corners (e.g., $\psi|_{\Gamma_{1}}(C)$ ,$\psi|_{\Gamma_{2}}(C)$ and $\psi|_{\Gamma_{3}}(C)$ in Fig.~\ref{fig:explain}) but on different patches (type B and C). 

\section{Numerical Results}
To validate the proposed explicit continuity enforcement scheme for EFIE discretized by CBIE method, we consider plane wave scattering from a PEC sphere, a cube, and a complex CAD model, and we compare the results obtained from our proposed scheme with those produced by solving the EFIE with the original CBIE proposed in~\cite{hu2021chebyshev} without continuity enforcement, as well as with the MFIE. We also show the computed radar cross section of a dipole antenna excited by a current sheet placed at the gap of the antenna and compare against that computed using the commercial Ansys HFSS solver and the MFIE, demonstrating that the proposed scheme can handle practical scenarios. 
\subsection{PEC sphere}
We first consider the canonical problem of scattering from a PEC sphere with diameter $D=2\lambda$ in free space illuminated by a plane wave excitation given by $\mathbf{E}^{\text{inc}} = \mathbf{e}_{x}e^{ik_{0}z}$. The surface is partitioned into 6 non-overlapping patches and GMRES is used to obtain the solution density. We also solved this problem by using the original CBIE formulation without continuity enforcement with and without taking into account the strongly singular contour integral, which may not always evaluate to 0 in this scenario due to the continuity not being enforced. In table~\ref{tab:cond}, the forward map (application of the integral operator to a prescribed density) convergence of the three approaches is shown. The number of GMRES iterations required is also shown along with the solution error when solving the matrix system. The GMRES residual tolerance was set to $10^{-6}$ from $N=10$ to $N=16$, and $10^{-9}$ for $N=18$ and $N=20$. It can be seen that both solvers without continuity enforcement take significantly more iterations regardless of whether the contour integral is included or not. Moreover, the original CBIE formulation without including the contour integral results in worse solved density accuracy for the same number of points per patch compared both other cases, while the forward map error for all three cases is similar, indicating the necessity for properly dealing with the continuity issue for Nystr\"om methods even for smooth geometries. The proposed explicitly enforced continuity approach requires significantly fewer iterations to reach convergence than both of the other methods due to reduction of the feasible solution search space by removing all spurious solutions that do not satisfy the normal continuity conditions.

\begin{table*}[b]
    \centering
    \caption{Convergence vs. points per patch per dimension (N) for EFIE with/without continuity enforcement }
    \label{tab:cond}
    \begin{tabular}{| c | c | c | c | c | c | c | c | c | c |}
        \hline
         & \multicolumn{3}{c|}{ ``Closed" CBIE} & \multicolumn{3}{c|}{``Open" CBIE}&
         \multicolumn{3}{c|}{``Open" + Contour integral}\\
         \hline
        $N$ & Fwd. Map Err.  & \# GMRES Iter. & Sol. Err.  & Fwd. Map Err.  & \# GMRES Iter. & Sol. Err. & Fwd. Map Err.  & \# GMRES Iter. & Sol. Err.\\
        \hline   \hline
        10 & 3.40e-03  & 112  & 2.72e-02 & 9.81e-03  & 654 & 2.09e+00 &  7.51e-04  & 472 & 9.36e-03  \\
        12 & 1.54e-04  & 82  & 1.07e-03 & 6.46e-04  & 756 & 8.75e-02 & 4.53e-05 & 605  &6.99e-04 \\
        14 & 1.06e-05  & 24 & 9.71e-06 & 3.00e-05  & 706 & 2.17e-03 & 1.88e-06  & 755 & 6.86e-05 \\
        16 & 9.75e-07 & 18 & 1.46e-06 & 3.27e-06  & 428 & 3.79e-04 & 1.58e-07 & 911 & 4.81e-05  \\
        18 & 1.29e-07 & 118 & 4.05e-07 & 4.18e-07 & 1534 & 1.74e-04 & 2.01e-08 & 1424 & 4.67e-07 \\
        20 & 1.90e-08  & 58  & 6.04e-08 & 6.34e-08  & 1710  & 3.83e-05 &5.31e-09 & 1658  & 1.31e-07\\
        \hline
    \end{tabular}
\end{table*}

\subsection{PEC cube}
\begin{figure*}
  \centering\includegraphics[width=0.8\textwidth,height=4.8cm]{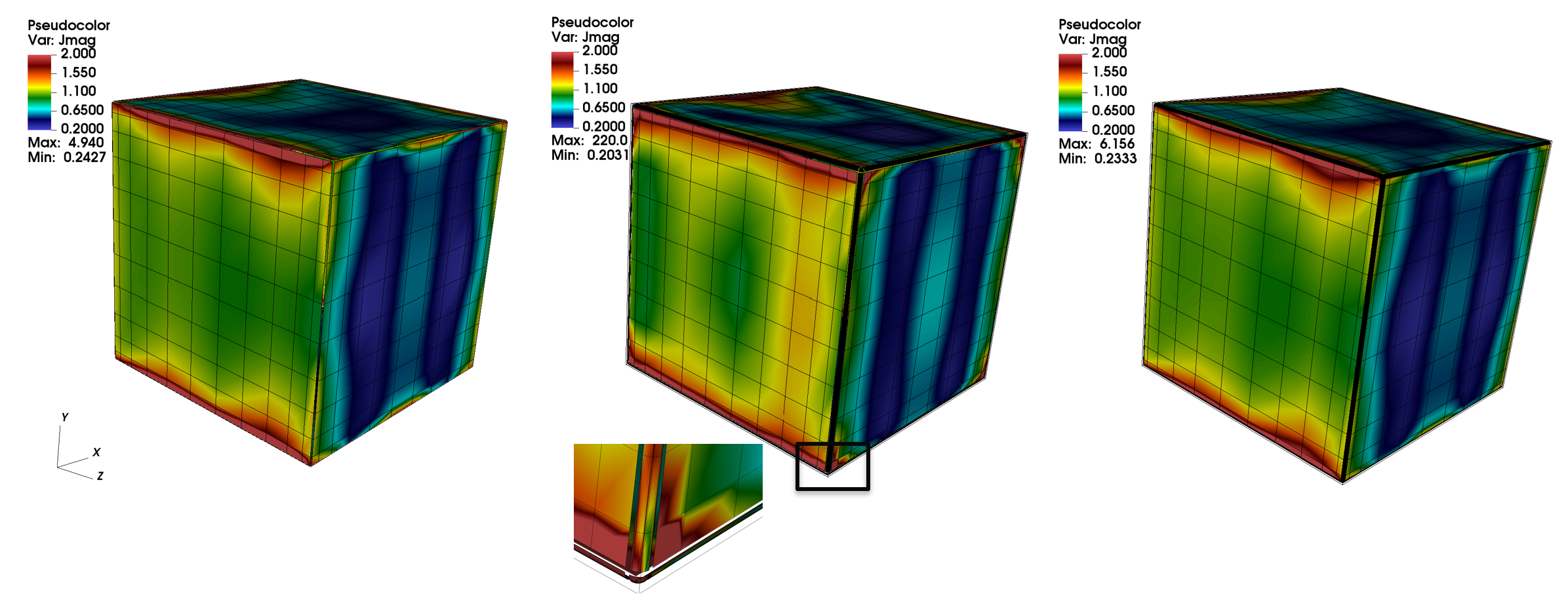}
  \caption{(Left) Current magnitude solved with MFIE. (Middle) Current magnitude solved with EFIE without continuity enforcement. (Right) Current magnitude solved with EFIE with explicit continuity enforcement.}
  \label{fig:cube}
\end{figure*}

The second example we present is the scattering from a PEC cube with 1$\lambda$ edge length. To avoid edge singular behavior, the edges and corners are slightly rounded by using quarter-cylinder patches with $0.01\lambda$ radius and 1/8th sphere patches with the same radius respectively. The same excitation is used as before. Three different solution approaches were compared: MFIE using the original CBIE with only interior quadrature points, the EFIE using the same approach without continuity enforcement of the densities, and EFIE using the proposed new CBIE that includes points on edges and corners to explicitly enforce the continuity. Fig.~\ref{fig:cube} shows the magnitude of the current density solved with each formulation. It can be seen that the results obtained with the EFIE with the proposed continuity enforcement agree well with the MFIE but differ significantly with the EFIE without continuity enforcement. In particular, hot spots at the corners are observed for EFIE without continuity enforcement and the values at those points appear to blow up, which indicates that the solution near the corner is incorrect. 

\subsection{PEC CAD model}
\begin{figure}[t]
\centering
\subfloat[][]{
\includegraphics[width=40mm]{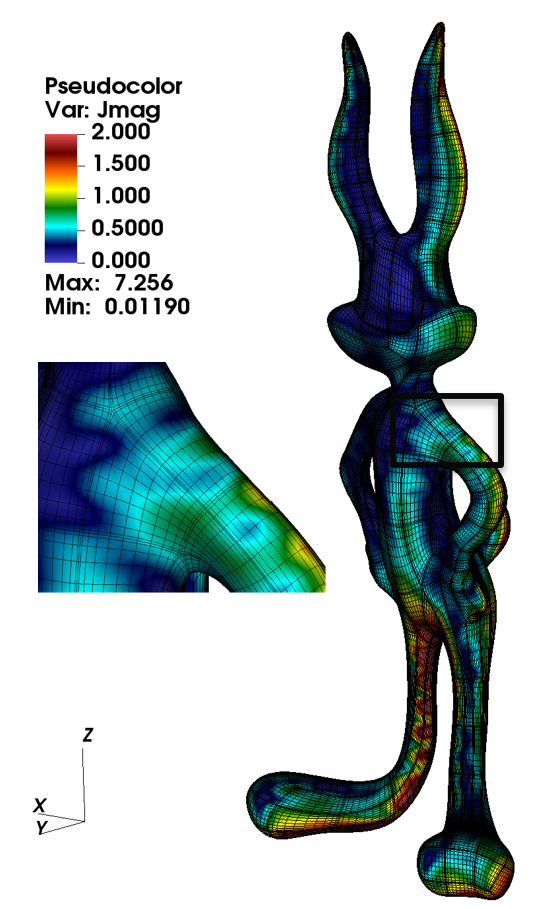}
} 
\subfloat[][]{
\includegraphics[width=38mm]{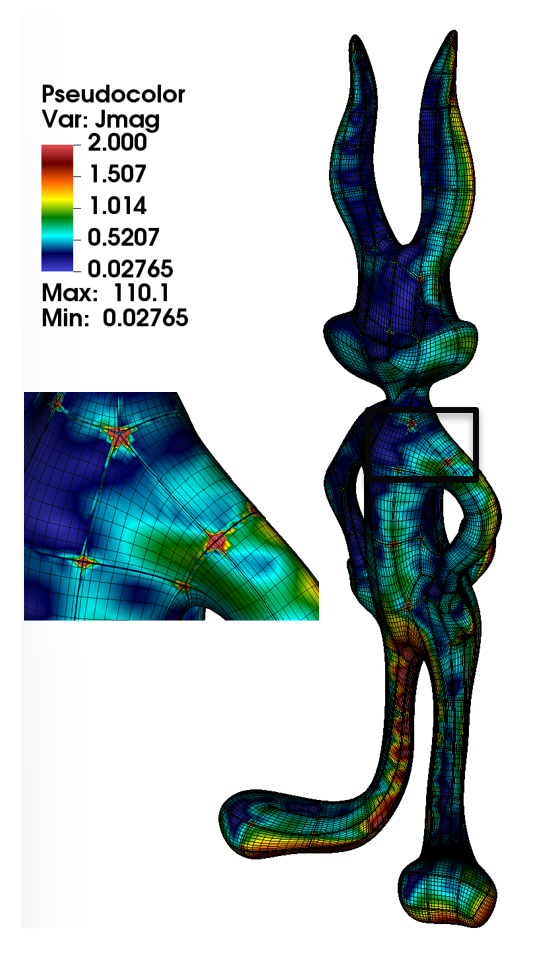}
} 

\caption{Surface electric current density induced on $16\lambda$ tall PEC CAD humanoid bunny model by incident plane wave solved with (a) EFIE with explicit continuity enforcement (b) EFIE without continuity enforcement}
\label{fig:CAD_bugs}
\end{figure}

To demonstrate that our explicit continuity enforcement scheme can also readily handle complex geometries, we consider scattering from a 16 wavelength tall NURBS humanoid bunny model produced by CAD software. The excitation is the same plane wave as used in the previous two examples. Fig.~\ref{fig:CAD_bugs} shows the current density solved by EFIE with and without explicit continuity enforcement respectively. It can be seen that the density values blow up at the intersections of different patches if the current continuity is not enforced despite these patch boundaries not corresponding to any geometric edges, whereas this issue is not present in the solution produced by the new proposed continuity enforcement scheme.

\subsection{PEC dipole antenna}
\begin{figure}[ht]
\centering
\subfloat[][]{
\includegraphics[width=42mm]{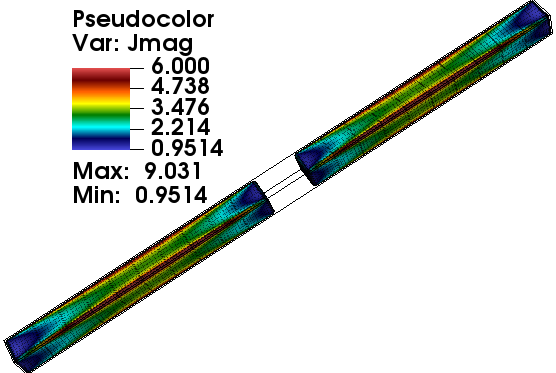}
}
\subfloat[][]{
\includegraphics[width=45mm]{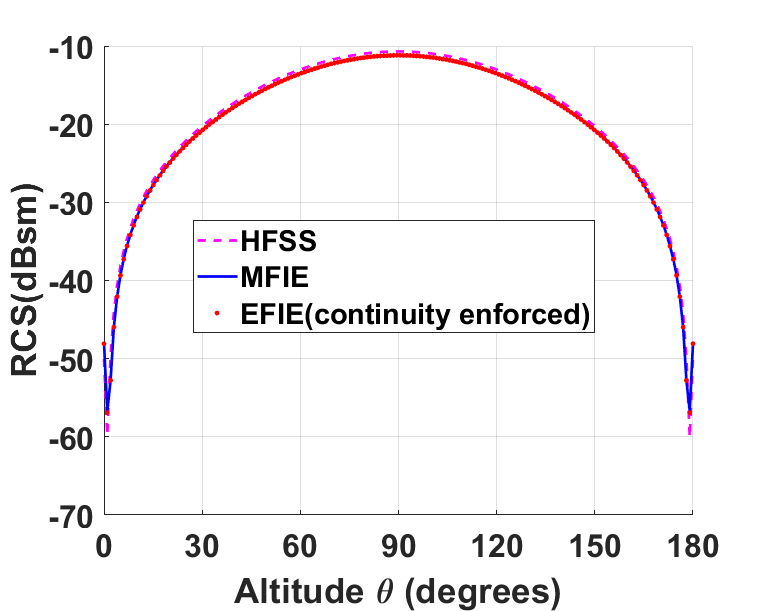}
} 
\caption{(a) Surface electric current density induced on the surface of a dipole structure by an incident plane wave. The total length of the dipole is $0.54\lambda$. (b) RCS at $\phi=0^{\circ}$  corresponding to plane wave scattering obtained with MFIE, EFIE with continuity enforced and HFSS.}
\label{fig:CAD_dipole}
\end{figure}

A dipole antenna structure is analyzed with the proposed scheme. This structure consists of two cuboid PEC arms, each with $0.025\lambda \times 0.025\lambda$ rectangular cross section and $0.25\lambda$ length. The gap between the arms is $0.04\lambda$. The corners and edges are slightly rounded in the same way as the previous example using a $0.0025\lambda$ rounding radius. The antenna is excited with a plane wave propagating in $+x$ direction and polarized in $-z$ direction. Fig.~\ref{fig:CAD_dipole} (a) and (b) shows the current density distribution on the dipole surface and the far-field pattern at the E-plane ($\phi= 0^{\circ}$) computed with Ansys HFSS, the EFIE with CBIE and continuity enforcement, and the MFIE with original CBIE respectively. It can be seen that the EFIE results and MFIE results are completely overlapped with each other and they both match very closely with the HFSS results.
\section{Conclusion}
In this letter, we extended the work in~\cite{hu2021chebyshev} to be applicable for solving the EFIE for PEC object scattering, which required explicit continuity enforcement of the densities across patch boundaries. We presented a simple and effective method to achieve this continuity by switching to a Clenshaw-Curtis quadrature which includes the endpoints for integration. Several numerical examples were presented and compared against the MFIE formulation and Ansys HFSS, showing excellent agreement in all cases and verifying the effectiveness of the proposed scheme.
\bibliographystyle{IEEEtran}
\bibliography{IEEEabrv,bibliography} 

\end{document}